\documentclass[11pt]{article}
\usepackage{graphicx,amssymb,subfigure}

\newcommand{\GLip}{{\rm \sc LipGr}}
\newcommand{\Median}{{\rm Median}}

\begin{document}

\title{Correction for ``Connect The Dots: How Many Random Points Can A Regular Curve Pass
Through?''}
\author{E. Arias-Castro, D.~L. Donoho, X. Huo and C.~A. Tovey}
\date{{\it Adv. in Appl. Probab.} {\bf 37}, no. 3 (2005), 571--603}
\maketitle


We are grateful to Dr. Eitan Bachmat of Ben-Gurion University, whose insightful comments lead to the present correction.
 
\begin{itemize}
\item
In Section 5.2, we state that
$$\Median(N_n) \sim 2 \sqrt{n},$$
while the correct statement should be
$$\Median(N_n) \sim \sqrt{2 n}.$$

\item
In Section 5.3, Figure 11 misstates the computational
experiment that was performed.
The experiment used the class of
Lipschitz graphs with Lipschitz constant 2.
However, the paper erroneously
stated that the experiment used
the class of  Lipschitz graphs with
constant 1.
In other words, in the experiment we computed $N_n(\GLip_2)/\sqrt{n}$ instead of $N_n(\GLip_1)/\sqrt{n}$ as stated in both the caption of Figure 11 and the comments in Section 5.3.\\ 

In the following figure we computed $N_n(\GLip_1)/\sqrt{n}$, and we immediately notice that, indeed, it converges to $\sqrt{2}$.

\begin{center}
\includegraphics[height=5in]{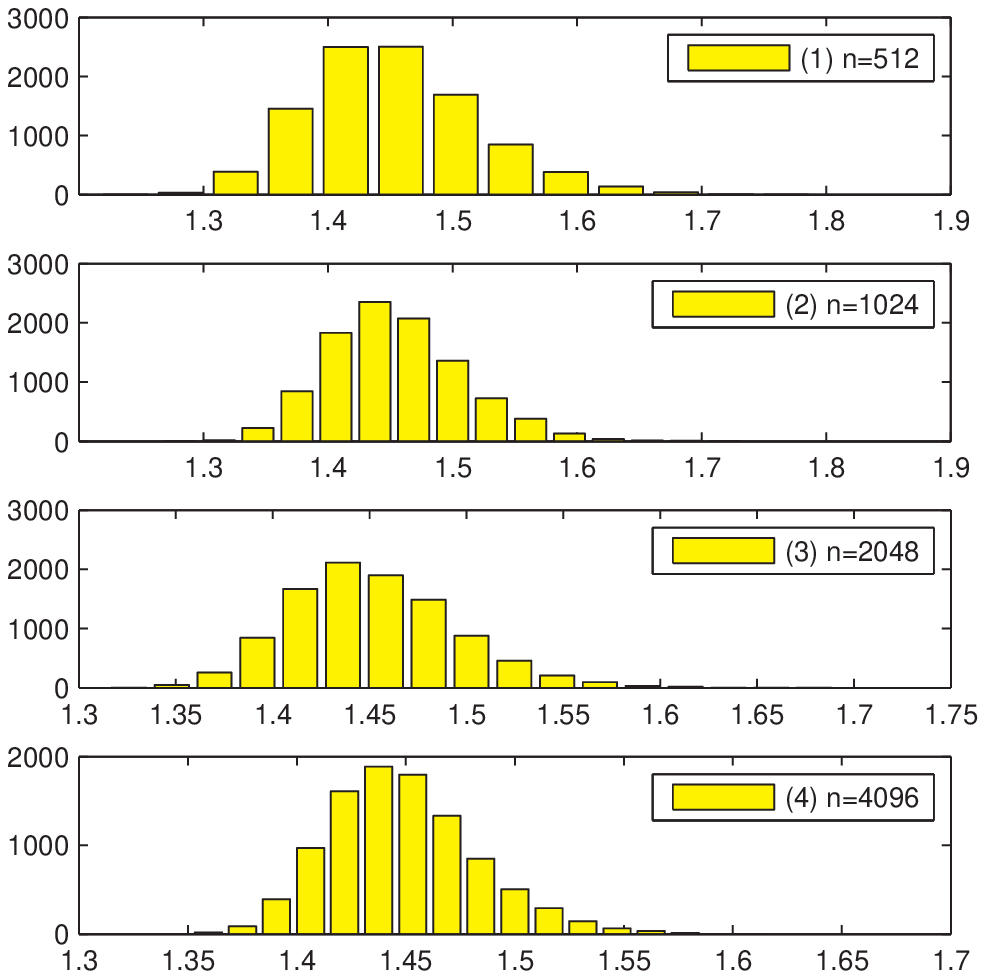}
\end{center}

\end{itemize}




\end{document}